\documentstyle{amsppt}
\magnification1200
\NoBlackBoxes
\def\phi{\varphi}
\def\RE{\text{\rm Re }}

\def\lap{{\Cal L}}
\def\hbar{\overline{h}}

\def\hchi{\hat \chi}
\def\hsigma{\hat \sigma}
\def\hbar{\overline{h}}
\pageheight{9 true in}
\pagewidth{6.5 true in}

\topmatter
\title
The number of unsieved integers up to $x$
\endtitle
\author Andrew Granville and K. Soundararajan
\endauthor
\rightheadtext{The number of unsieved integers up to $x$}
\address{D{\'e}partment  de
Math{\'e}matiques et Statistique, Universit{\'e}
de Montr{\'e}al, CP 6128 succ
Centre-Ville, Montr{\'e}al, QC  H3C 3J7, Canada}\endaddress
\email{andrew{\@}dms.umontreal.ca}
\endemail
\address{Department of Mathematics, University of Michigan, Ann Arbor,
Michigan 48109, USA} \endaddress \email{ksound{\@}umich.edu} \endemail
\thanks{Le premier auteur est partiellement soutenu par une bourse
du Conseil  de recherches en
sciences naturelles et en g\' enie du Canada. The
second  author is partially supported by the National Science Foundation and the
American Institute of Mathematics (AIM).}
\endthanks

\abstract
Typically, one expects that there are around $x\prod_{p\not\in P,\ p\leq x}
(1-1/p)$ integers up to $x$,
all of whose prime factors come from the set $P$. Of
course for some choices of $P$ one may get rather
more integers, and for some choices
of $P$ one may get rather less.
Hall [4] showed that one never gets more than
$e^\gamma+o(1)$ times the expected amount
(where $\gamma$ is the Euler-Mascheroni
constant), which was improved slightly by Hildebrand [5].
 Hildebrand [6] also showed that for a
given value of $\prod_{p\not\in P,\ p\leq x} (1-1/p)$, the smallest count that you get
(asymptotically) is when $P$ consists of all the primes up to a given point. In this
paper we shall improve Hildebrand's upper bound, obtaining a result close to optimal,
and also give a substantially shorter proof of Hildebrand's lower bound. As part of
the proof we give an improved Lipschitz-type bound for such counts.
\endabstract

\endtopmatter
\document

\head 1. Introduction \endhead

\noindent Let $f$ denote a multiplicative function with $0\le f(n)\le 1$
for all positive integers $n$.  Define
$$
\Theta(f,x):= \prod_{p\leq x} \left( 1-\frac 1p\right)
\left( 1+\frac {f(p)}p+\frac{f(p^2)}{p^2}+\dots \right).
$$
Given a real number $w\ge 1$ in this paper we are concerned with
the problem of determining
$$
g(w):=\liminf_{x\to \infty} \frac{1}{x}\sum_{n\le x} f(n), \qquad
\text{and} \qquad G(w):=\limsup_{x\to \infty} \frac{1}{x}
\sum_{n\le x}f(n), \tag{1.1}
$$
where both limits are taken over the class of multiplicative
functions $f$ with $\Theta(f,x) = 1/w+o(1)$.

P. Erd{\H o}s and I. Ruzsa [1] showed that $g(w)>0$ for all $w$.
Consider the function $f$ with $f(p^k) = 1$ for $p\le x^{1/w}$ and
$f(p^k)=0$ for $x^{1/w}\le p \le x$.  Then one has
$\Theta(f;x)=1/w+o(1)$ and further $\sum_{n\le x} f(n)=
\psi(x,x^{1/w}) $, the number of integers below $x$ having no
prime factors above $x^{1/w}$. It is well known that for any fixed
$w$ we have
$$
\psi(x,x^{1/w}) = x \rho(w) \Big(1 +O\Big(\frac{w}{\log
x}\Big)\Big), \tag{1.2}
$$
where $\rho(w)$ is the Dickman--de Bruijn function, defined by
$\rho(w) =1$ for $0 \le w \le 1$, and $w\rho^{\prime}(w) =
-\rho(w-1)$ for all $w \ge 1$. This example shows that $g(w) \le
\rho(w)$ and A. Hildebrand [6] established that in fact $g(w)=
\rho(w)$.  Since $\rho(w) = w^{-w+o(w)}$ note that $g(w)$ decays
very rapidly as $w$ increases.

Regarding $G(w)$, R. Hall [4] established that $G(w) \le e^{\gamma}/w$ and Hildebrand
[5] improved this slightly by showing that $G(w)\le \frac{1}{w}\int_0^w \rho(t) dt$.
Since $\int_0^{\infty} \rho(t)dt =e^{\gamma}$ this does mark an improvement over
Hall's result, but the difference from $e^{\gamma}/w$ is $\frac{1}{w}\int_{w}^{\infty}
\rho(t) dt = w^{-w+o(w)}$ which is very small.   In this paper we shall prove that
$G(w) = e^{\gamma}/w -1/w^{2+o(1)}$, but it remains to determine $G(w)$ more
precisely. We shall also give a shorter proof of Hildebrand's result that
$g(w)=\rho(w)$.

\proclaim{Theorem 1} For all $w\ge 1$ we have that
$$
G(w) \ge \max_{w\ge \Delta \ge 0} \Big(\rho(w+\Delta) +
\int_{0}^{\Delta} \frac{\rho(t)}{w+\Delta -t} dt\Big). \tag{1.3}
$$
When $w$ is large, the maximum is attained for $\Delta \sim \log
w/\log \log w$, and yields
$$
G(w) \ge \frac{e^{\gamma}}w - \frac{(e^{\gamma}+o(1))\log w} {w^2
\log \log w}.
$$
\endproclaim

\proclaim{Theorem 2} For all large $w$ we have
$$
G(w)\le \frac {e^\gamma}w - \frac 1{w^2 \exp (c (\log w)^{2/3}
(\log\log w)^{1/3}) }
$$
for a positive constant $c$.
\endproclaim

We also give an explicit upper bound for $G(w)$ valid for all $w$.

\proclaim{Theorem 3} For $1\le w$ we have that $G(w) \le 1-\log w
+(\log w)^2/2$ and equality holds here for $1\le w\le 3/2$.  For
$w\ge 1$ put $\Lambda(w):= \frac 12 (w+1/w) + \frac{\log
w}{2}(w-1/w)$.  Then $G(w)\le \Lambda(w) \log
(1+e^{\gamma}/(w\Lambda(w)))$.
\endproclaim

The first bound in Theorem 3 is better than the second for $w\le 3.21\ldots$, when the
second bound takes over. Note that the second bound in Theorem 3 equals
$e^\gamma/w-(e^{2\gamma}+o(1))/w^3\log w$, only a little weaker than the bound in
Theorem 2, while being totally explicit.

 In the range $1\le w\le 3/2$ we may check that the right side
of (1.3) equals $1-\log w +(\log w)^2/2 =G(w)$.
 Perhaps it is true that $G(w)$ is
given by the right side of (1.3) for all $w$.

We end this section by giving a simple construction that proves
Theorem 1.

\demo{Proof of Theorem 1} Let $y$ be large and consider the
completely multiplicative function $f$ defined by $f(p)=0$ for
$p\in [y,y^w]$ and $f(p)=1$ for all other primes $p$.  Put
$x=y^{w+\Delta}$ where $0\le \Delta \le w$ and note that
$\Theta(f,x)=\prod_{y\leq p\leq y^w} (1-1/p) \sim 1/w$.  An
integer $n\le x$ with $f(n)=1$ has at most one prime factor
between $y^w$ and $x$, and all its other prime factors are below
$y$.  Hence
$$
\sum_{n\le x} f(n) = \psi(x,y) + \sum_{y^w\le p\le x} \psi(x/p,y),
$$
and using (1.2) and the prime number theorem this is
$$
\sim x\rho(w+\Delta) + x\sum_{y^w \le p\le x} \frac{1}{p}
\rho\Big(w+\Delta -\frac{\log p}{\log y}\Big) \sim x \Big(
\rho(w+\Delta) + \int_0^{\Delta} \frac{\rho(t)}{w+\Delta -t}
dt\Big),
$$
which gives the lower bound (1.3) for $G(w)$.  For large $w$ we
see that
$$
\rho(w+\Delta) +\int_{0}^{\Delta} \frac{\rho(t)}{w+\Delta-t} dt =
\frac{1}{w+\Delta} \int_{0}^{\Delta} \rho(t)dt + \int_0^{\Delta}
\frac{t\rho(t)}{(w+\Delta)(w+\Delta-t)}dt + \rho(w+\Delta)
$$
and since $\int_0^{\infty} t\rho(t)dt <\infty$ and
$\int_0^{\Delta} \rho(t) dt =e^{\gamma} -
\Delta^{-(1+o(1))\Delta}$ the above is
$$
\frac{1}{w+\Delta} (e^{\gamma} -\Delta^{-(1+o(1))\Delta})
+O\Big(\frac{1}{w^2}\Big).
$$
The quantity above attains a maximum for $\Delta =(1+o(1))\log
w/\log \log w$, completing the proof of Theorem 1.
\enddemo

We noted above that $G(w)=1-\log w+(\log w)^2/2$ for $1\leq w\leq 1.5$ (with the
maximum attained in (1.3) at $\Delta=w$).  Next we record the bounds obtained for
$1.5\leq w\leq 2$ (though here the maximum is attained with $\Delta$ a little smaller
than $w$).

\medskip
\centerline{\vbox{\offinterlineskip \hrule \halign{& \vrule \ # & \ \vrule \,
\vrule\ # & \ \vrule \ #& \ \vrule \ #& \ \vrule \ #&\ \vrule \ #&\ \vrule \ # &\strut
\ \vrule #\cr $w$ \hfil &\hfil$1.5$ &\hfil$1.6$ &\hfil$1.7$ &\hfil$1.8$ &\hfil$1.9$
&\hfil$2.0$ &\cr \noalign{\hrule}\cr \noalign{\hrule} $G(w)\geq$\hfil
&\hfil.676735&\hfil.640255&\hfil.608806&\hfil.581685&\hfil.557392&\hfil.535905&\cr
\noalign{\hrule}\cr $G(w)\leq$ \hfil
&\hfil.676736&\hfil.640449&\hfil.610155&\hfil.584960&\hfil .564135&\hfil.547080&\cr
\noalign{\hrule}\cr} \hrule}}

\smallskip

\centerline{\sl The upper and lower bounds for $G(w)$ given by Theorems 1 and 3.}
\medskip

\head 2.  Reformulation in terms of integral equations \endhead

\noindent  E. Wirsing [8] observed that questions on mean-values of
multiplicative functions can be reformulated in terms of solutions
to a certain integral equation.  We formalized this connection
precisely in our paper [2] and we now recapitulate the salient
details.  We will prove our results by establishing the
corresponding statements for solutions to integral equations.

The following class of integral equations is relevant to the study of multiplicative
functions $f$ with $|f(n)|\leq 1$ for all $n$: Let $\chi$ be a measurable function
with $\chi(t) =1$ for $t\le 1$ and $|\chi(t)|\le 1$ for all $t\ge 1$. Let
$\sigma(u)=1$ for $u\le 1$ and for $u>1$ we define $\sigma$ to be the solution to
$$
u\sigma(u) = \int_0^u \chi(t) \sigma(u-t) dt. \tag{2.1}
$$
In [2] we showed that there is a unique continuous solution
$\sigma(u)$ to (2.1) and that $|\sigma(u)| \le 1$ for all $u$.  In
fact $\sigma(u)$ is given by
$$
\sigma(u) = 1+ \sum_{j=1}^{\infty} \frac{(-1)^j}{j!} I_j(u;\chi),
\tag{2.2a}
$$
where
$$
I_j(u;\chi) =  \int\Sb t_1,\ldots, t_j\ge 1\\ t_1+\ldots +t_j \le
u\endSb \frac{1-\chi(t_1)}{t_1} \cdots \frac{1-\chi(t_j)}{t_j}
dt_1\cdots dt_j . \tag{2.2b}
$$

The connection between multiplicative functions and the integral
equation (2.1) is given by the following result which is
Proposition 1 in [2].

\proclaim{Proposition 2.1} Let $f$ be a multiplicative function
with $|f(n)|\le 1$ for all $n$, and $f(n)=1$ for $n\le y$. Let
$\vartheta(x) =\sum_{p\le x} \log p$ and define
$$
\chi(u) = \chi_f(u) = \frac{1}{\vartheta(y^u) } \sum_{p\le y^u}
f(p)\log p.
$$
Then $\chi(t)$ is a measurable function with $|\chi(t)|\le 1$ for
all $t$ and $\chi(t)=1$ for $t\le 1$. Let $\sigma(u)$ be the
corresponding unique solution to (2.1). Then
$$
\frac{1}{y^u} \sum_{n\le y^u} f(n) =\sigma(u)
+O\biggl(\frac{u}{\log y}+\frac{1}{y^u}\biggr).
$$
\endproclaim

For our problems on non-negative multiplicative functions we will
restrict our attention to integral equations where $\chi(t)$ only
takes values in $[0,1]$.  The corresponding solution $\sigma(u)$
to (2.1) then also takes only values in $[0,1]$.  We also define
$$
E(u)= E_{\chi}(u) := \exp\Big(\int_0^u \frac{1-\chi(t)}{t} dt
\Big). \tag{2.3}
$$
Notice that $\Theta(f,y^u)\sim E(u)$ when $\chi$ is defined as in Proposition 2.1.
Analogously to $g(w)$ and $G(w)$ we may define
$$
{\tilde g}(w) = \liminf\Sb u, \chi\\ E_{\chi}(u)=w\endSb
\sigma(u), \qquad\text{and} \qquad {\tilde G}(w)= \limsup\Sb u,
\chi \\ E_{\chi}(u)=w \endSb \sigma(u),
$$
where the limits are taken over all $\chi$ with $\chi(t)=1$ for
$t\le 1$ and $\chi(t)\in [0,1]$ for all $t$, and over all $u\ge 1$
with $E_\chi(u)=w$.  We shall show that these quantities are in
fact equal to $g(w)$ and $G(w)$ respectively.  Something similar
was stated (but not very precisely) by Hildebrand in his
discussion paper [7].

\proclaim{Theorem 2.2} We have $g(w)={\tilde g}(w)$ and $G(w)
={\tilde G}(w)$.
\endproclaim

To prove Theorem 2.2 we need to know how small primes affect the
mean-values of multiplicative functions, so that we can
remove their effect and be in a position to use Proposition 2.1.  We
also require a converse to Proposition 2.1 which allows us
to go from integral equations to multiplicative functions.
Such results were established in [2] and we now quote them in
our context.  Proposition 4.4 of [2]
(with $\phi=\pi/2$ there) gives the following Lemma.

\proclaim{Lemma 2.3}  Let $f$ be a multiplicative function with
$0\le f(n)\le 1$ for all $n$.  Let $1\ge \epsilon \ge \log 2/\log x$
and take $g$ to be the completely multiplicative function with $g(p)=1$ if
$p\leq x^\epsilon$, and $g(p)=f(p)$ otherwise.
Then
$$
\frac 1x \sum_{n\le x} f(n)= \Theta(f,x^{\epsilon})
\frac 1x \sum_{m\le x} g(m) + O( \epsilon^{\frac 14 - \frac{1}{2\pi}}).
$$
\endproclaim

Next, Proposition 1 (Converse) in [2] gives the following converse to
Proposition 2.1.

\proclaim{Proposition 2.4} Let $\chi$ be a given measurable function with $\chi(t)=1$
for $t\le 1$ and $\chi(t)\in [0,1]$ for all $t \ge 1$ and let $\sigma$ denote the
corresponding solution to {\rm (2.1)}.  Given $\epsilon>0$ and $u\ge 1$ there exist
arbitrarily large $y$ and a multiplicative function $f$ with $f(n)=1$ for $n\le y$ and
$0\le f(n) \le 1$ for all $n$ and with
$$
\Big|\chi(t) - \frac{1}{\vartheta(y^t)}\sum_{p\le y^t} f(p)\log p\Big|
\le \epsilon
\qquad \text{for almost all } 0\le t\le u.
$$
Further, for all $1\le t\le u$
$$
\sigma(t)= \frac{1}{y^t} \sum_{n\le y^t} f(n) + O(u^{\epsilon} -1) +
O\Big(\frac{u}{\log y}\Big).
$$
\endproclaim

We will defer the proof of Theorem 2.2 to the
next section.  But let us note that combining Lemma 2.3 with Proposition 2.1
gives
$$
g(w) \ge \min_{w\ge v\ge 1} \frac{1}{v} {\tilde g}\Big(\frac wv\Big),
\qquad \text{and} \qquad
G(w) \le \max_{w\ge v\ge 1} \frac{1}{v} {\tilde G}
\Big(\frac wv\Big). \tag{2.4a}
$$
Also from Proposition 2.4 we get that
$$
g(w) \le {\tilde g}(w), \qquad \text{and} \qquad
{\tilde G}(w) \le G(w).
\tag{2.4b}
$$

We end this section by recording two facts which will be
useful in our later work.  Firstly when $\chi(t) \in [0,1]$
one obtains inclusion-exclusion inequalities from (2.2a,b): namely,
for all even integers $n$ we have (see Proposition 3.6 of [2])
$$
\sum_{j=0}^{n} \frac{(-1)^j}{j!} I_j(u;\chi) \ge \sigma(u)\ge
\sum_{j=0}^{n+1} \frac{(-1)^j}{j!} I_j(u;\chi). \tag{2.5}
$$
Secondly from (2.2a,b) and a little combinatorics we obtain the following Lemma (see
Lemma 3.4 of [2]):

\proclaim{Lemma 2.5}  Let $\chi$ and ${\hat \chi}$ be two measurable functions with
$\chi(t)=\hchi(t)=1$ for $0\le t\le 1$ and $|\chi(t)|$, $|\hchi(t)|\le 1$ for all $t$.
Let $\sigma$ and $\hsigma$ be the corresponding solutions to {\rm (2.1)}. Then
$\hsigma(u)$ equals
$$
\sigma(u) + \sum_{j=1}^{\infty} \frac{1}{j!} \int\Sb t_1,\ldots ,
t_j \ge 1\\
t_1+\ldots +t_j \le u \endSb \frac{\hchi(t_1)-\chi(t_1)}{t_1} \ldots
\frac{\hchi(t_j)-\chi(t_j)}{t_j} \sigma(u-t_1-\ldots-t_j)dt_1\ldots dt_j.
$$
\endproclaim

\head 3. Upper bounds for $G(w)$ and Lipschitz estimates \endhead

\noindent For a measurable function $g:[0,\infty) \to {\Bbb C}$ we will
denote the Laplace transform of $g$ by
$\Cal L (g,s):=\int_0^\infty g(t)e^{-st} dt$. If $g$ is integrable and grows
sub-exponentially (that is, for every $\epsilon >0$, $|g(t)| \ll_{\epsilon}
e^{\epsilon t}$ almost everywhere)
then the Laplace transform  is well defined for all complex numbers $s$
with $\RE(s) >0$.  Integrating term by term in  (2.2a,b) we
see that
$$
{\Cal L}(\sigma,s)= \frac{1}{s}
\exp\Big(-\lap\Big(\frac{1-\chi(v)}{v},s\Big)\Big). \tag{3.1}
$$

Suppose now that
$\chi(t) =1$ for $t\le 1$ and $\chi(t)\in [0,1]$ for all $t$ and we are
given $u\ge 1$.  Define $\hchi(t)=\chi(t)$ for $t\le u$ and $\hchi(t)=0$
for $t>u$.  If $\sigma$ and $\hsigma$ are the corresponding solutions
to (2.1) then note that $\sigma(v)=\hsigma(v)$ for $v\le u$ and
that $E_{\chi}(v)=E_{\hchi}(v)$ for $v\le u$.  Now
$$
\sigma(u) = \hsigma(u) = \frac{1}{u}\int_0^u \hsigma(t)\hchi(u-t)dt
\le \frac{1}{u}\int_0^u \hsigma(t) = \frac 1u \int_0^\infty \hsigma(t) dt
- \frac 1u \int_u^{\infty} \hsigma(t)dt.
$$
Further
$$
\align
\lap\Big(\frac{1-\hchi(t)}{t},s\Big)-\log E(u)
&= \int_0^\infty \left(\frac{1-\hchi(t)}{t} \right) e^{-st} dt
- \int_0^u  \frac{1-\chi(t)}t \ dt \\
&=
\int_0^u \left( \frac{1-\chi(t)}{t} \right) (e^{-st}-1) dt  +  \int_u^\infty
\frac{e^{-st}}{t}  dt \\ &= -\gamma -\log (su) + O( u |s|) , \\
\endalign
$$
for small $s$, since
$\gamma = \int_0^{1} \frac{1-e^{-t}}{t} dt - \int_1^{\infty}
\frac{e^{-t}}{t} dt$.  Hence
$$
 \frac{1}{u}\int_0^{\infty} \hsigma(t) dt = \frac{1}{u} \lim_{y\to 0}
\lap(\hsigma,y) = \lim_{y\to 0} \frac{1}{yu}
\exp\Big(-\lap\Big(\frac{1-\hchi(t)}{t},y\Big)\Big)
 = \frac{e^{\gamma}}{E(u)},
$$
and so we have
$$
\sigma(u) \le \frac{e^{\gamma}}{E(u)} - \frac 1u \int_u^{\infty}
\hsigma(t) dt. \tag{3.2}
$$
We  use (3.2) in the proofs of Theorems 2 and 3, since it allows us
 to give an upper bound for $\sigma(u)$ by determining a "smoothed lower bound"
for $\hsigma$. Our plan for proving a bound on this integral is to bound how much
 $\hsigma(t)$ changes as $t$ gets bigger than $u$, via Lipschitz-type estimates.

For general  complex $\chi$ with $|\chi|\le 1$, and $\sigma$ satisfying (2.1)
we might expect to have a Lipschitz estimate of the form
$$
\Big||\sigma(u)|-|\sigma(v)|\Big| \ll \Big(\frac{u-v}{u}\Big)^{\kappa} \left( 1 + \log
\frac{u}{u-v} \right) \ \text{whenever} \ 1\le v\le u \tag{3.3}
$$
for certain values of $\kappa >0$; and indeed we established (3.3) in [3] for $\kappa
=1-2/\pi$. Any increase in the value of $\kappa$ allows stronger consequences, and we
believe that $\kappa=1$ in (3.3) is probably valid. Note that no exponent $>1$ is
possible since $|\rho(1+\delta)-\rho(1)|=\log (1+\delta) \sim \delta$ for
$0\le \delta= \le 1$. 
We are able to improve ``$1-2/\pi$'' to ``$1-1/\pi$'' in the special
case that $\chi(t)\in [0,1]$ for all $t$.

\proclaim{Theorem 4}  Let $\chi$ be a measurable function with $\chi(t) =1$ for
$t\le 1$ and $\chi(t)\in [0,1]$ for $t>1$, and let $\sigma$ denote the corresponding
solution to {\rm (2.1)}.  Then
$$
|\sigma(u)-\sigma(v)|
\ll \Big(\frac{u-v}{u}\Big)^{1-\frac{1}{\pi}}
 \Big( 1 + \log \frac{u}{u-v}\Big) \ \text{whenever} \
1\le v\le u.
$$
\endproclaim

Theorem 4 follows immediately from the stronger but more complicated Proposition 4.2
below, and the  fact that $|\sigma(u)-\sigma(v)|\leq  \frac{3(u-v)}u$ whenever $v\leq
u(1-1/E(u))$.
This is trivial for $v\leq 2u/3$, whereas for larger $v$ in the range,
 we obtain
$$
|\sigma(u)-\sigma(v)|\leq \frac{e^\gamma}{E(v)} \leq \frac{ue^\gamma}{vE(u)} \leq
\frac{3(u-v)}u,
$$
using Hall's result that $\sigma(u)\leq e^\gamma/E(u)$.

Using (3.3) in (3.2) leads to the bound ${\tilde G}(w)
\le e^{\gamma}/{w} - C_\kappa/(w^{1+1/\kappa}\log w)$ for some
positive constant $C_\kappa$.  
Thus if (3.3) holds with $\kappa=1$ then we would be able to 
deduce that $G(w)=e^\gamma/w-(\log w)^{O(1)}/w^2$ by Theorem 1.

In order to prove Theorem 3 we give the following
explicit Lipschitz estimate (see also Proposition 4.1 of [2]).

\proclaim{Proposition 3.1} Let $\chi$ be a measurable function with $\chi(t)=1$ for
$t\le 1$ and $\chi(t) \in [0,1]$ for all $t$, and let $\sigma(u)$ denote the
corresponding solution to {\rm (2.1)}.  Then for all $u\ge 1$ and $1\ge \delta >0$ we
have
$$
\log(1+\delta) \Big(\frac{E(u)-1/E(u)}{2} + \log E(u)
\frac{E(u)+1/E(u)}{2}\Big)
\ge
\sigma(u(1+\delta)) -\sigma(u),
$$
and
$$
\sigma(u(1+\delta))-\sigma(u) \ge - \log(1+\delta)
\Big(\frac{E(u)+1/E(u)}{2} + \log E(u) \frac{E(u)-1/E(u)}{2}\Big).
$$
\endproclaim
\demo{Proof} We shall only prove the lower bound, the proof
of the upper bound is similar.  From (2.2a,b) we see that
$$
\sigma(u(1+\delta)) -\sigma(u) \ge -
\sum\Sb j=1\\ j \text{ odd }\endSb^{\infty} \frac{1}{j!} \left(
I_j(u(1+\delta);\chi)-I_j(u;\chi)\right).
$$
By symmetry we see that $I_j(u(1+\delta);\chi)-I_j(u;\chi)$ equals
$$
j \int\Sb
t_1,\ldots,t_{j-1}\ge 1 \endSb \frac{1-\chi(t_1)}{t_1} \cdots
\frac{1-\chi(t_{j-1})}{t_{j-1}} \int\Sb
\max(t_1,\ldots,t_{j-1},u-t_1-\ldots-t_{j-1})\le t_j \\
t_j\le u(1+\delta)-t_1-\ldots-t_{j-1}\endSb \frac{1-\chi(t_j)}{t_j}
dt_1\cdots dt_j.
$$
The integral over $t_j$ is
$$
\le \log \frac{u/j+u\delta}{u/j} = \log (1+j\delta) \le j\log (1+\delta),
$$
since $\max(t_1,\ldots,t_{j-1},u-t_1-\ldots-t_{j-1}) \ge u/j$.  Further
since $\delta <1$ we have $t_1$, $\ldots$, $t_{j-1} \le u$ and
so these integrals contribute $\le (\log E(u))^{j-1}$.  Thus we
have
$$
\sigma(u(1+\delta)) -\sigma(u) \ge -
\sum\Sb j=1\\ j \text{ odd }\endSb^{\infty} \frac{1}{j!} j^2\log (1+\delta)
(\log E(u))^{j-1},
$$
and the result follows easily.

\enddemo

\demo{Proof of Theorem 2.2} Fix $w\geq v\geq 1$. Suppose $\chi(t)=1$ for $t\le 1$ and
$\chi(t)\in [0,1]$ for all $t$ and let $\sigma(u)$ denote the corresponding solution
to (2.1) (we will think of $\chi$ as giving the optimal function for either ${\tilde
g}(w/v)$ or ${\tilde G}(w/v)$). Let $U\ge 1$ be a parameter which we will let tend to
infinity. Put $\chi_1(t)=\chi(t/U)$ and note that the corresponding solution to (2.1)
is $\sigma_1(u) =\sigma(u/U)$. Define $\chi_2(t)=0$ for $1\le t\le v$ and $\chi_2(t)
=\chi_1(t)$ for all other $t$, and let $\sigma_2(u)$ denote the corresponding solution
to (2.1).  By Lemma 2.5 we see that for $U\ge v$
$$
\sigma_2(uU) = \sigma_1(uU) +\sum_{j=1}^{\infty} \frac{(-1)^j}{j!}
\int\Sb v\ge t_1,\ldots, t_j \ge 1 \\ t_1+\ldots+t_j \le uU\endSb
\frac{1}{t_1} \frac{1}{t_2} \ldots \frac{1}{t_j}
\sigma_1(uU-t_1-\ldots -t_j) dt_1 \cdots dt_j.
$$
By Proposition 3.1 we know that
$$
\sigma_1(uU-t_1-\ldots -t_j)
= \sigma_1(uU) + O\Big(\min\Big(1,E_{\chi}(u)
\log E_{\chi}(u) \frac{jv}{uU}\Big)\Big).
$$
Using this above we see easily that for large $U$ with $u,v,w$ fixed we have
$\sigma_2(uU) \sim \sigma_1(uU)/v=\sigma(u)/v$ and note further that $E_{\chi_2}(uU) =
vE_{\chi_1}(uU)= vE_{\chi}(u)$.

This scaling argument shows that for $1\le v\le w$ we have
${\tilde g}(w/v) \ge v{\tilde g}(w)$ and that
${\tilde G}(w/v) \le v {\tilde G}(w)$.  Using these inequalities
in (2.4a) we deduce that $g(w)\ge {\tilde g}(w)$ and that $G(w)\le
{\tilde G}(w)$ and combining this with (2.4b) we obtain Theorem 2.2.

\enddemo

Now that Theorem 2.2 has been established, to prove Theorem 3
it suffices to establish the analogous bounds for ${\tilde G}(w)$
and we establish these next.

\demo{Proof of Theorem 3}  Using the inclusion-exclusion upper bound (2.5)
with $n=2$ we see that $\sigma(u) \le 1-\log E(u) + (\log E(u))^2/2$.
It follows that $G(w)={\tilde G}(w) \le 1-\log w+ (\log w)^2/2$.
If $w\le 3/2$ then consider $\chi(t)=0$ for $1\le t\le w$ and $\chi(t)=1$
for all other $t$.  Then we see that the corresponding solution $\sigma(u)$
satisfies $\sigma(u)=1-\log w +(\log w)^2/2$ for $3\ge u\ge 2w$.  Thus
${\tilde G}(w)=1-\log w +(\log w)^2/2$ for $1\le w\le 3/2$.

We now establish the second bound of the Theorem.  As noted in the
introduction the second bound is worse than the first for $w\le 3.21$ and
so we may suppose that $w\ge 2$.  With $\hchi, \hsigma$ as above,
note that $\hsigma(t)\ge 0$ for all $t$, and
$$
\hsigma(u(1+\delta))\ge \hsigma(u) - \Lambda(E(u)) \log (1+\delta) \
\text{for} \ 0\le \delta \le 1
$$
by Proposition 3.1. If $E(u)\ge 2$ then $\Lambda(E(u))\ge 7/4 > 1/\log 2$ so that
$\exp(\sigma(u)/\Lambda(E(u)))-1 <1$.  Hence we obtain that
$$
\align
\frac 1u\int_{u}^{\infty} \hsigma(t)dt &\ge
\int_0^{\exp(\sigma(u)/\Lambda(E(u)))-1} (\sigma(u)-\Lambda(E(u))
\log(1+\delta)) d\delta\\
&= -\sigma(u) +\Lambda(E(u)) \Big(\exp\Big(\frac{\sigma(u)}{\Lambda(E(u))}
\Big)-1\Big),\\
\endalign
$$
and inserting this into (3.2) we get the Theorem.
\enddemo


\head 4. An improved upper bound: Proof of Theorem 2 \endhead

\noindent Our proof of Theorem 2 is also based on (3.2) and obtaining
lower bounds for $\frac 1u\int_u^{\infty} \hsigma(t) dt$. However Theorem 4
is not quite strong enough to obtain this conclusion and so, in this
section, we develop a hybrid Lipschitz estimate which for our problem
is almost as good as
(3.3) with $\kappa =1$. We begin with the following Proposition (compare Lemma 2.2 and
Proposition 3.3 of [3]).

\proclaim{Proposition 4.1}  Let $\chi$ be a measurable function with $\chi(t)=1$ for
$t\le 1$ and $\chi(t)$ in the unit disc for all $t$. Let $\sigma$ be the corresponding
solution to {\rm (2.1)}. Let $1\le v\le u$ be given real numbers, and put $\delta =
u-v$.  Define
$$
F:= \max_{y\in {\Bbb R}}  \  \exp\left( \gamma-\int_0^u \text{\rm Re} \left(
\frac{1-\chi(t) e^{-ity}}{t}\right) dt\right) \ |1-e^{-iy\delta}| .
$$
Then
$$
\align |\sigma(u)-\sigma(v)| &\le \frac{\delta}{u} \log \frac{eu}{\delta}+ F + F
\int_0^{2/(uF)} \frac{1-e^{-2xu}}{x} dx \\
& \le \frac{\delta}{u} \log \frac{eu}{\delta}+ F\log \frac {e^3} F .\\
\endalign
$$
\endproclaim

\demo{Proof} As in the proof of Theorem 3 take $\hchi(t)=\chi(t)$
for $t\le u$ and $\hchi(t)=0$ for $t>u$, and let $\hsigma$ be
the corresponding solution to (2.1).  Set $\sigma(t)=\hsigma(t)=0$
for $t<0$. Note that
$$
\align |u\sigma(u)-v\sigma(v)| &= |u\hsigma(u)-v\hsigma(v)|
= \Big|\int_0^u \chi(t) (\hsigma(u-t) - \hsigma(v-t)) dt\Big| \\
&\le \int_0^u |\hsigma(t) -\hsigma(t-\delta)| dt = \int_0^u 2t
|\hsigma(t)-\hsigma(t-\delta)| \Big(\int_0^{\infty} e^{-2xt} dx
\Big) dt\\
&\le 2 \int_0^{\infty} \int_0^u \{
|t\hsigma(t) - (t-\delta)\hsigma(t-\delta)|
+\delta|\hsigma(t-\delta)|\} e^{-2tx} dt dx\\
&\le \int_0^{\infty} I(x) dx + \int_0^{\infty}
 \int_{\delta}^u 2\delta e^{-2tx} dt dx
= \delta \log \frac{u}{\delta}
+  \int_0^{\infty} I(x) dx,\\
\endalign
$$
where
$$
I(x) = \int_0^{u} 2|t\hsigma(t) - (t-\delta)\hsigma(t-\delta)| e^{-2tx} dt.
$$
As $|\sigma(u)-\sigma(v)| \le \frac{1}{u} ( |u\sigma(u)-v\sigma(v)| +
\delta|\sigma(v)|) \le \frac{\delta}{u} + \frac{1}{u} |u\sigma(u)- v\sigma(v)|$, it
follows that
$$
|\sigma(u)-\sigma(v)| \le \frac{\delta}{u}\log \frac{eu}{\delta} + \frac{1}{u}
\int_0^{\infty} I(x) dx. \tag{4.1}
$$

By Cauchy's inequality
$$
\align I(x)^2 &\le \Big( 4\int_0^{u} e^{-2tx} dt \Big) \Big( \int_0^{u}
|t\hsigma(t)-(t-\delta)\hsigma(t-\delta)|^2 e^{-2tx} dt\Big) \\
&\le 2 \Big(\frac{1-e^{-2xu}}{x} \Big) \Big(\int_0^{\infty}
|t\hsigma(t)-(t-\delta)\hsigma(t-\delta)|^2 e^{-2tx} dt\Big).
\\
\endalign
$$
By Plancherel's formula the second term above is
$$
=\frac{1}{2\pi } \int_{-\infty}^{\infty} |\lap(t\hsigma(t) -(t-\delta)
\hsigma(t-\delta), x+iy)|^2 dy
= \frac{1}{2\pi} \int_{-\infty}^{\infty} |\lap(t\hsigma(t),x+iy)|^2
|1-e^{-(x+iy)\delta}|^2 dy.
$$
From (2.1) we see that $\lap(t\hsigma(t),x+iy)
= \lap(\hsigma,x+iy)\lap(\hchi,x+iy)$
and so the above equals
$$
\frac{1}{2\pi }\int_{-\infty}^{\infty}
|\lap(\hsigma,x+iy) \lap (\hchi,x+iy)|^2
|1-e^{-(x+iy)\delta}|^2 dy
\le F(x)^2 \cdot \frac{1}{2\pi} \int_{-\infty}^{\infty}
|\lap(\hchi,x+iy)|^2 dy
$$
where
$$
F(x):=\max_{y\in {\Bbb R}}|1-e^{-(x+iy)\delta}| |\lap(\hsigma,x+iy)|.
$$
Now,  using Plancherel's formula again,
$$
\frac{1}{2\pi}
\int_{-\infty}^{\infty} |\lap(\hchi,x+iy)|^2 dy = \int_0^{\infty}
|\hchi(t)|^2 e^{-2tx} dt \leq  \int_0^{u} e^{-2tx} dt
= \frac{1-e^{-2xu}}{2x},
$$
and so
$$
I(x) \le  \frac{1-e^{-2xu}}{x} F(x). \tag{4.2}
$$

We now demonstrate that $F(x)$ is a decreasing function of $x$.
Suppose that $\beta>0$ is real, and
recall that the Fourier transform of $k(z):=e^{-\beta|z|}$ is
${\hat k}(\xi) = \int_{-\infty}^{\infty}
e^{-\beta|z|-i\xi z}dz=\frac{2\beta}{\beta^2+\xi^2}$.
Hence $e^{-\beta z} = k(z)= k(-z)= \frac{1}{\pi}\int_{-\infty}^{\infty}
\frac{\beta}{\beta^2 +\xi^2}e^{-i\xi z}dz$ by Fourier inversion for
$z>0$.  It follows
that for $\delta+t>0$ we have
$$
(1-e^{-\delta(x+\beta+iy)})e^{-t(x+\beta+iy)}
= \frac{1}{\pi} \int_{-\infty}^{\infty}
\frac{\beta}{\beta^2+\xi^2} e^{-t(x+iy+i\xi)}(1-e^{-\delta(x+iy+i\xi)})d\xi.
$$
Multiplying both sides by $\hsigma(t)$, and
integrating $t$ from $0$ to $\infty$, we
deduce that
$$
\align (1-e^{-\delta(x+\beta+iy)})\lap(\hsigma,x+\beta+&iy) = \frac{1}{\pi}
\int_{-\infty}^{\infty} \frac{\beta}{\beta^2+\xi^2} \lap(\hsigma,x+iy+i\xi)
(1-e^{-\delta(x+iy+i\xi)}) d\xi \\
&\le \Big( \max_{y\in {\Bbb R}}
|(1-e^{-\delta(x+iy)})
\lap(\hsigma,x+iy)|\Big) \frac{1}{\pi}\int_{-\infty}^{\infty}
\frac{\beta}{\beta^2+\xi^2} d\xi,\\
\endalign
$$
and so $F(x+\beta)\leq F(x)$ as  claimed. Therefore
$F(x) \le   \lim_{x\to 0^+} F(x)$.

Now if $s=x+iy$ with $x>0$ then
$$
\align \lap\biggl(\frac{1-\chi(v)}{v},s\biggr)&= \int_0^\infty \left(
\frac{1-\chi(v)e^{-ivy}}v \right) e^{-vx} dv
+ \int_0^\infty \frac{e^{-vs}-e^{-vx}}v
dv \cr
&= \int_0^\infty \left( \frac{1-\chi(v)e^{-ivy}}v \right) e^{-vx} dv  +
\log(x/s),\cr
\endalign
$$
so that
$$
\lap(\sigma,s)
 =\frac{1}{x}
\exp\biggl(-\int_0^\infty \left( \frac{1-\chi(v)e^{-ivy}}v
\right) e^{-vx} dv \biggr).
$$
Using this for $\hsigma$ we have
$$
|\lap(\hsigma,x+iy)| = \frac{1}{x}\exp\Big( -\int_{u}^{\infty}
\frac{e^{-tx}}{t} dt -
\int_0^u \text{Re}
\Big( \frac{1-\chi(t)e^{-ity}}{t} \Big) e^{-tx} dt\Big).
$$
For $x\ll 1/u$ we get
$$
\int_{u}^{\infty} \frac{e^{-tx}}{t} dt
= \int_{ux}^{\infty} \frac{e^{-t}}{t} dt  =
\int_1^{\infty} \frac{e^{-t}}{t} dt
+ \int_{ux}^{1} \frac{e^{-t}-1}{t} dt + \log
\frac{1}{ux} = -\gamma +\log \frac{1}{ux} + O(ux),
$$
since
$\gamma = \int_0^{1} \frac{1-e^{-t}}{t} dt
- \int_1^{\infty} \frac{e^{-t}}{t}dt$,
so that
$$
|\lap(\hsigma,x+iy)|
=e^{\gamma} u\exp\Big(-\int_0^u \text{Re}
\left(\frac{1-\chi(t)e^{-ity}}{t} \right)dt + O(ux) \Big) .
$$
Note that this is $\ll_u 1$, so that the maximum of $|1-e^{-(x+iy)\delta}|
|\lap(\hsigma,x+iy)|$
cannot occur with $\| y\delta/2\pi \|\to 0$ as $x\to 0^+$ (here
$\| t\|$ denotes the distance from the nearest integer to $t$), else
$F(x)\ll_u x+\| y\delta/2\pi \| \to
0$ as $x\to 0^+$, implying that $F(x)=0$ which is ridiculous.
Thus the maximum occurs
with $\| y\delta/2\pi \|\gg 1$ as $x\to 0^+$ so that $1-e^{-(x+iy)\delta} =
1-e^{-iy\delta}+O(x\delta)=(1-e^{-iy\delta}) \{ 1+O(x\delta)\}$, so that
$$
|1-e^{-(x+iy)\delta}| |\lap(\hsigma,x+iy)| = u |1-e^{-iy\delta}|
 \exp\left( \gamma-\int_0^u \text{\rm Re} \left(
\frac{1-\chi(t) e^{-ity}}{t}\right) dt +O(ux) \right) .
$$
Therefore
$F(x) \le   uF\{ 1+O(ux)\}$ for sufficiently small $x$; and so $F(x) \le
uF$. Also $F(x) \le 2\max_{y\in {\Bbb R}}|\lap(\hsigma,x+iy)| \le 2/x$.
Therefore, by
(4.2), we get that
$$
I(x)\le \cases
\frac{1-e^{-2xu}}{x} uF &\text{if } x\le 2/uF\\
\frac{2}{x^2} &\text{if } x>2/uF,\\
\endcases
$$
which when inserted in (4.1) yields the first estimate in the Proposition.

Now if $F\leq 1$ then
$$
\int_0^{2/(uF)} \frac{1-e^{-2xu}}{x} dx
\leq \int_0^{2/u} \frac{1-e^{-2xu}}{x} dx +
\int_{2/u}^{2/(uF)} \frac{1}{x} dx \leq 2+\log (1/F) ,
$$
and so we deduce the second estimate of Proposition 4.1.
If $F>1$ this holds trivially
since $|\sigma(u)-\sigma(v)|\leq 2$.

\enddemo

As an application of this Proposition, we establish the following
strange-looking Lipschitz estimate in the case that
$\chi(t) \in [0,1]$ for all $t\ge 1$.

\proclaim{Proposition 4.2}  Let $\chi$ be a measurable function with $\chi(t) =1$ for
$t\le 1$ and $\chi(t)\in [0,1]$ for $t>1$, and let $\sigma$ denote the corresponding
solution to {\rm (2.1)}. Let $1\le v\le u$ be given  and write $E(u)=(u/(u-v))^\theta$
for $\theta>0$.  Then
$$
|\sigma(u)-\sigma(v)| \ll \Big(\frac{u-v}{u}\Big)^{ \min \{ 1,1-\frac{1}{\pi} \sin(\pi
\theta)\}} \Big( 1 + \log \frac{u}{u-v}\Big).
$$
\endproclaim

\demo{Proof} Let $\delta=u-v$ and $A=\int_0^u \frac{1-\chi(t)}{t} dt=\log E(u)$. We
will show that
$$
 \exp\Big(-\int_0^u \frac{1-\chi(t) \cos(ty)}{t} dt\Big) \min(1,\delta y)
\ll \Big(\frac{\delta}{u}\Big)^{ \min
\{ 1,1-\frac{1}{\pi} \sin(\frac{\pi A}{\log
(u/\delta)})\}} , \tag{4.3}
$$
for all positive $y$. The result then follows from Proposition 4.1 since $F\ll$ Left
side of (4.3).

If $y\le e/u$ then the left side of (4.3) is
$\leq e\delta/u$ and the result follows.
Henceforth we may suppose that $y>e/u$.
Since $\cos(x)=1+O(x^2)$, we get that
$\int_0^{1/y} \frac{1-\chi(t)\cos(ty)}{t} dt
= \int_0^{1/y} \frac{1-\chi(t)}{t} dt
+O(1)$. Thus if we let
$z:=\int_{1/y}^{u} \frac{1-\chi(t)}{t} dt$ then
$$
\align \int_{0}^{u} \frac{1-\chi(t)\cos(ty)}{t} dt
&= A-z+O(1) + \int_{1/y}^{u}
\frac{1-\chi(t)\cos(ty)}{t} dt
\\
&= A-z+O(1) + \int_{1/y}^{u} \frac{1-\cos (ty)}{t} dt + \int_{1/y}^{u}
\frac{1-\chi(t)}{t} \cos(ty) dt\\
&=A-z+\log (uy) + O(1) + \int_{1}^{uy} \frac{1-\chi(t/y)}{t} \cos(t) dt,\\
\endalign
$$
by making a change of variables, and since (integrating by parts)
$$
\int_{1/y}^{u} \frac{\cos(ty)}{t} dt = \frac{\sin(ty)}{yt}\Big|_{1/y}^{u}
+ \int_{1/y}^{u} \frac{\sin(ty)}{yt^2} dt = O(1).
$$
By periodicity
$$
\int_1^{uy}\frac{1-\chi(t/y)}{t} \cos(t) dt
=\int_{0}^{\pi} G(\theta) \cos\theta \
d\theta, \ \ \text{where} \
G(\theta):= \sum\Sb t \pm \theta \in 2\pi\Bbb Z \\ 1\le
t\le uy \endSb \frac{1-\chi(t/y)}{t}
$$
and the sum over $t$ above is over real values of $t$ in the
range $[1,uy]$ such that
$t \pm \theta$ is an integer multiple of $2\pi$.  Note that
$$
\align 0\le G(\theta)
\le \frac{1}{\pi} \log &(uy)+O(1) \ \text{ for all } \ \theta ,
\\
\text{and} \ \
\int_0^{\pi} G(\theta) d\theta &= \int_{1/y}^{u} \frac{1-\chi(t)}{t} dt
= z.  \\
\endalign
$$
Consider the problem of minimizing
$\int_0^{\pi} G(\theta) \cos \theta d\theta$ over
all functions $G$ satisfying these two constraints.
Since $\cos \theta$ decreases from
$1$ to $-1$ in the range $[0,\pi]$, we see that this is achieved by taking
$G(\theta)=0$ for $\theta \in [0,\pi-\theta_0]$, and $G(\theta)
=\frac{1}{\pi} \log
(uy) +O(1)$ for $\theta \in [\pi-\theta_0,\pi]$, where
$\theta_0$ satisfies $\theta_0
(\frac{1}{\pi } \log (uy) +O(1)) =z$. We conclude that
$$
\align \int_{0}^{\pi } G(\theta)  \cos \theta  d\theta
&\ge \int_{\pi-\theta_0}^{\pi}
\cos \theta \Big(\frac{1}{\pi}\log (uy)+O(1)\Big) d\theta
= -\frac{1}{\pi} \log (uy) \sin \theta_0 + O(1)\\
&= -\frac{1}{\pi} \log (uy) \sin \Big(\frac{\pi z}{\log(uy) +O(1)}\Big) +O(1)
\\
&=  -\frac{1}{\pi} \log (uy) \sin \Big(\frac{\pi z}{\log(uy)}\Big) +O(1),
\\
\endalign
$$
since $0\leq z\leq \log (uy)$. Therefore
$$
\int_{0}^{u} \frac{1-\chi(t)\cos(ty)}{t} dt
\ge A-z+\log (uy) \Big(1- \frac{1}{\pi}
\sin \Big(\frac{\pi z}{\log(uy)}\Big)\Big) +O(1). \tag{4.4}
$$
In the domain $0\leq z\leq \log (uy)$, the right side of (4.4)
is a non-increasing
function of $z$, so that it is
greater than the value with $z$ replaced by $\log(uy)$,
that is, it is $>A+O(1)$. Therefore the left side of (4.3)
is $\ll e^{-A} \min(1,\delta y)$,
which is $\leq  \delta/u$ if $A\geq \log (uy)$, as required.
If $A<\log (uy)$ then the
right side of (4.4) is greater than the value with
$z$ replaced by $A$, which is $\log(uy) - \frac{\log (uy)}{\pi}
\sin(\pi A/\log(uy)) +O(1)$, so that the left side of (4.3)
is
$$
\ll \frac{\min(1,\delta y)}{uy} (uy)^{\frac 1\pi
\sin(\frac{\pi A}{\log (uy)})}.
$$
This function is maximized when $y=1/\delta$ in the range
$\log (uy) \ge A$, at which
point it yields the right side of (4.3),
completing the proof.
\enddemo

\demo{Proof of Theorem 2}  Let $\alpha=E(u)=e^A$.
We may assume that $\alpha$ is
large, and that $\sigma(u) \ge 1/\alpha$, else our result follows
trivially.  Let
$v=(1+e^{-\lambda})u$ for some parameter $\lambda >A$,
and select $\hchi(t)=\chi(t)$ for $t\le u$ and $\hchi(t)=0$ for $t>u$,
as earlier.
Using Proposition 4.2 we deduce that there is a constant $C$ such that
$$
|\hsigma(u)-\hsigma(v)| \le C (1+\lambda) \exp\Big( -\lambda +
\frac{\lambda}{\pi}\sin\Big(\frac{\pi A}{\lambda}\Big)\Big). \tag{4.5}
$$
If $\lambda \ge 2A$, then this is $\le C(1+\lambda)
\exp(-\lambda(1-1/\pi))$ which is
easily verified to be $\le 1/(2\alpha)$ if $\alpha$ is sufficiently large. If
$A<\lambda \le 2A$, then the right side of (4.5)
is $\le 2C(1+A) \exp(-\lambda
+\frac{\lambda}{\pi} \sin(\frac{\pi A}{\lambda}))$,
which is a decreasing function of
$\lambda$ in our range. For
$\lambda=A+\xi$ where $\xi:=c A^{2/3} (\log A)^{1/3}$,
with $c>(6/\pi^2)^{1/3}$, this equals
$$
2C (1+A) \exp\Big(-A-\xi+\frac{A+\xi}{\pi}
\sin\Big(\frac{\pi A}{A+\xi}\Big) \Big) =
2C (1+A) \exp\Big(-A -\frac{\pi^2}{6} \frac{\xi^3}{A^2} + O\left(
\frac{\xi^4}{A^3}\right) \Big)
 \le
\frac{1}{2\alpha}.
$$
Thus we have proved that $|\hsigma(u)-\hsigma(v)|\le 1/(2\alpha)$
for all $\lambda
\geq A+\xi$, which implies that $\hsigma(v)\ge 1/(2\alpha)$
for $u\leq v\leq
u(1+e^{-A-\xi})$. Therefore
$$
\frac{1}{u} \int_{u}^{\infty} \hsigma(t) dt \ge \frac{1}{u}
\int_{u}^{u(1+e^{-A-\xi})}
\hsigma(v) dv \ge \frac{1}{u} \cdot ue^{-A-\xi} \cdot \frac{1}{2\alpha} >
\frac{1}{2\alpha^2 \exp(\xi)},
$$
which implies the theorem, by (3.2).

\enddemo

\head 5.  Determining $g(w)={\tilde g(w)}$: Preliminaries \endhead

\noindent In the remainder of the paper we will give an
alternative, substantially shorter, proof of Hildebrand's
result that $g(w)={\tilde g}(w)=\rho(w)$.  More precisely,
we will establish the following Theorem.

\proclaim{Theorem 5}  Let $\chi(t)=1$ for $t\le 1$ and $\chi(t)\in [0,1]$ for all
$t>1$, and let $\sigma(u)$ denote the corresponding solution to {\rm (2.1)}.  Then
$\sigma(u) \ge \rho(E(u))$ for all $u$.  Further if $1\le E(u)\le 2$ and
$\sigma(u)=\rho(E(u))$ then $E(u/2)=1$. If $E(u)\ge 2$ and $\sigma(u)=\rho(E(u))$ then
$E(u/E(u))=1$; that is, $\chi(t)=1$ for 
$t\leq u/E(u)$, and $\chi(t)=0$ for $u/E(u)\le t\le u$, except possibly on 
a set of measure $0$.
\endproclaim

If $1\le E(u)\le 2$ then using (2.5) with $n=0$ we see
that $\sigma(u)\ge 1- I_1(u;\chi) =1-\log E(u) = \rho(E(u)$.
Further (2.5) with $n=2$ gives that
$$
\sigma(u) \ge 1-I_1(u;\chi)+ \frac 12 I_2(u;\chi) - \frac 16 I_3(u;\chi)
\ge 1-\log E(u) + \frac 12 I_2(u;\chi) \Big(1-\frac{\log E(u)}{3}\Big)
$$
so that $\sigma(u) =\rho(E(u))$ if and only if $I_2(u;\chi)=0$, or
in other words $E(u/2)=1$.  This proves Theorem 5 in the
range $1\le E(u)\le 2$ and we assume below that $E(u)>2$.

Henceforth we let $u_0 :=u/E(u) <u_1 :=u(1-1/E(u))$.  We also
define
$$
B(u) = B_{\chi}(u) = \int_0^u \chi(v) dv.
$$
We note a simple principle that we shall use repeatedly.

\proclaim{Lemma 5.1}
Let $b\ge a$ be real numbers.  Let $f:[a,b] \to [0,1]$ and
$g:[a,b]\to {\Bbb R}$ be measurable functions,
such that $g$ is non-decreasing in
$[a,b]$, with $A:=\int_a^b f(t) dt$. Then
$$
\int_a^{a+A} g(t) \le \int_{a}^{b} f(t) g(t) dt \le \int_{b-A}^b g(t) dt.
$$
\endproclaim
\demo{Proof}  To prove the lower bound note that
$$
\int_{a}^{b-A} f(t) g(t) dt
\leq g(b-A) \int_{a}^{b-A} f(t) dt = g(b-A) \int_{b-A}^{b}
(1-f(t))dt \leq \int_{b-A}^{b} g(t)(1-f(t))dt ,
$$
and the result follows. The upper bound can be proved analogously.
\enddemo

\proclaim{Lemma 5.2}  For all $0 \le t\le y$
$$
y \frac{E(t)}{E(y)} - t \le B(y)-B(t) \le y - t\frac{E(y)}{E(t)}.
$$
Written differently
$$
E(t) \le \frac{E(y)}{y} (t+B(y)-B(t)), \qquad \text{and} \qquad E(y) \le
\frac{E(t)}{t} (y-B(y)+B(t)).
$$
\endproclaim
\demo{Proof} Note that
$$
\frac{E(y)}{E(t)}
= \frac{y}{t} \exp\Big( -\int_t^y  \frac{\chi(v)}{v} dv\Big).
$$
Applying Lemma 5.1 (with $f(v)=\chi(v)$ and $g(v)=-1/v$) we deduce that
$$
-\log \frac{t+B(y)-B(t)}{t} \le - \int_t^y \frac{\chi(v)}{v} dv \le -\log
\frac{y}{y-B(y)+B(t)},
$$
and the Lemma follows.
\enddemo

We note that
$$
B(y)\ge \frac{y}{E(y)} , \quad \text{and} \quad \frac{E(t)}{t} \ge \frac{E(y)}{y} \ \
\text{for} \ 0\le t \le y , \tag{5.1}
$$
which is a  particular case of Lemma 5.2.





Our proof of Theorem 5 splits into two cases which we
handle by different methods.  The first case, which
we treat in section 6, is when either
$E(u)$ is small ($\le 2.6$) or if $E(u_0) \ge E(u)-1$ is large.
The other case concerns $E(u)\ge 2.6$ and $E(u_0) < E(u)-1$
which is handled in Section 7.

\head 6.  The case $2<E(u)\le 2.6$, or $E(u_0) \ge E(u)-1$
\endhead

\proclaim{Proposition 6.1}  If $E(u)>2$ and $E(u_0) \ge E(u)-1$ then
$\sigma(u) > \rho(E(u))$.
\endproclaim
\demo{Proof}  Define ${\hat \chi}(t)=\chi(t)$ for $t\le u_0$ and ${\hat \chi}(t)=1$
for $t>u_0$ and let $\hsigma$ denote the solution to the corresponding integral
equation.  By Lemma 2.5 we have
$$
\align
\sigma(u)-\hsigma(u)
&= \sum_{j=1}^{\infty} \frac{(-1)^j}{j!} \int\Sb t_1, \ldots ,t_j\ge u_0\\
t_1+\ldots+t_j \le u\endSb \frac{1-\chi(t_1)}{t_1} \ldots
\frac{1-\chi(t_j)}{t_j} \hsigma(u-t_1-\ldots -t_j) dt_1\ldots dt_j \\
&\ge - \sum_{j \text{ odd}} \frac{1}{j!} \Big(\int_{u_0}^{u}
\frac{1-\chi(t)}{t}dt\Big)^{j} =  - \frac{1}{2}
\Big(\frac{E(u)}{E(u_0)}-\frac{E(u_0)}{E(u)}\Big) \tag{6.1}\\
\endalign
$$

Let $2\le n$ denote the largest even integer below $E(u)$. In the integral defining
$I_j(u,\hchi)$ the integrand can be non-zero only if each $t_i\leq u_0$, so that
$0\leq I_j(u,\hchi)\leq (\log E(u_0))^j$ for all $j$.  Also we have $t_1+\dots
+t_j\leq ju_0\leq nu_0\leq u_0E(u)=u$ if $j\leq n$, implying that $I_j(u,\hchi)=(\log
E(u_0))^j$. Therefore by the inclusion-exclusion inequality (2.5) we see that
$$
 \hsigma(u) \ge 1+\sum_{j=1}^{n+1} \frac{(-1)^j}{j!}
 I_j(u,\hchi) \geq \sum_{j=0}^{n+1} \frac{(-1)^j}{j!} (\log E(u_0))^j .\tag{6.2a}
$$
 Further note that
$$
\align
\frac{I_{n+2}(u,\hchi)}{(n+2)!}  -\frac{I_{n+3}(u,\hchi)}{(n+3)!}
&\ge \frac{I_{n+2}(u,\hchi)}{(n+2)!} \Big(1-\frac{\log E(u_0)}{n+3}\Big)\\
&\ge \frac{(\log E(u/(n+2)))^{n+2}}{(n+2)!}
\Big(1-\frac{\log E(u_0)}{n+3}\Big) \\
&\ge \frac{(\log (E(u_0)E(u)/(n+2)))^{n+2}}{(n+2)!}
\Big(1-\frac{\log E(u_0)}{n+3}\Big),\\
\endalign
$$
since $E(u/(n+2)) \ge E(u_0) \frac{u}{(n+2)u_0} = E(u_0)E(u)/(n+2)$ by
(5.1).  Thus another lower bound furnished by (2.5) is
$$
\hsigma(u) \ge \sum_{j=0}^{n+1} \frac{(-1)^j}{j!} (\log E(u_0))^j
+ \frac{(\log (E(u_0)E(u)/(n+2)))^{n+2}}{(n+2)!}
\Big(1-\frac{\log E(u_0)}{n+3}\Big). \tag{6.2b}
$$

If $2\le E(u)\le 6$ then using (6.1) together
with (6.2b) for appropriate $n$ we checked that $\sigma(u) > \rho(E(u))$ if $E(u_0)\ge
E(u)-1$. If $6\le n\le E(u)\le n+2$ then the right side of (6.2a) is at least
$1/E(u_0) - (\log (n+2))^{n+2}/(n+2)!$ and combining this with (6.1) we get that for
$E(u_0) \ge E(u)-1$
$$
\sigma(u) \ge \frac{1}{2E(u)(E(u)-1)} - \frac{1}{(n+2)(n+1)}
\frac{(\log (n+2))^{n+2}}{n!}
\ge \frac{.014}{(n+2)(n+1)} > \rho(n)\geq \rho(E(u)),
$$
since $(\log (n+2))^{n+2}/n! \le (\log 8)^8/6! <.486$ for $n\ge 6$
(note that $.14/56> 2 \times 10^{-4}$ whereas
$\rho(6) \approx2 \times 10^{-5}$).

\enddemo

Henceforth we may assume that $E(u_0) \le E(u)-1$.  We complete
this section by giving a proof of Theorem 5 for the range
$2< E(u) \le 2.6$.

\proclaim{Proposition 6.2} If $2< E(u) \le 2.6$ then
$\sigma(u) \ge \rho(E(u))$ and equality holds only when $E(u_0)=1$.
\endproclaim
\demo{Proof} By Proposition 6.1
we may assume that $e^{\xi}:=E(u/3) \le E(u_0) \le E(u)-1$
so that  $\xi \le\log (E(u)-1)$.
Since one of $t_1$, $t_2$ or $t_3$ (in the definition of $I_3$) must
be less than $u/3$, we see easily that $I_3(u)
\le 3 (\int_{1}^{u/3}\frac{1-\chi(v)}{v}dv)
I_2(u) \le 3\xi I_2(u)$.  Thus using (2.5) with $n=2$ we get
$$
\sigma(u) \ge 1- I_1(u) + \frac{1}{2}(1-\xi)I_2(u) = 1-\log E(u) +
\frac{1}{2}(1-\xi)I_2(u). \tag{6.3}
$$

Now $I_1(u-v)/v$ is a non-increasing function so by Lemma 5.1 we obtain
$$
\align I_2(u) &\ge \int_{1}^{\frac u3} \frac{1-\chi(v)}{v} I_1(u-v) dv
+ \int_{\frac u3}^{u} \frac{1-\chi(v)}{v} I_1(u-v)dv\\
&\ge \int_{\frac u{3e^{\xi}}}^{\frac u3} I_1(u-v)\frac{dv}{v} +
\int_{\frac{ue^{\xi}}{E(u)}}^{u} I_1(u-v)\frac{dv}{v}.\\
\endalign
$$
Note that $I_1(t)
=\log E(t) \ge \log E(u) - \int_{t}^{u} dv/v = \log (E(u)t/u)$, and
also that $I_1(t) \ge \log E(u/3)=\xi$ if $t\ge u/3$.
Using these bounds above we get
$$
I_2(u)
\ge \int_{\frac u{3e^{\xi}}}^{\frac u3} \log \Big(\frac{E(u)}{u}(u-v)\Big)
\frac{dv}{v}
+ \int_{\frac{ue^{\xi}}{E(u)}}^{u_1} \log\Big(\frac{E(u)}{u}(u-v)\Big)
\frac{dv}{v} + \int_{u_1}^{2u/3} \xi \frac{dv}{v}.
$$
Let $\gamma(E(u))
= \int_{1}^{E(u)-1} \log (E(u)-t) \frac{dt}{t} = \int_{u_0}^{u_1}
\log (E(u)- vE(u)/u) \frac{dv}{v}$. We see that
$$
I_2(u) \ge
\gamma(E(u)) + \xi \log\Big(\frac{2}{3} \frac{E(u)}{E(u)-1}\Big) +
\int_{\frac u{3e^{\xi}}}^{\frac u3}
\log \Big(\frac{E(u)}{u}(u-v)\Big) \frac{dv}{v} -
\int_{\frac{u}{E(u)}}^{\frac{ue^{\xi}}{E(u)}}
\log \Big(\frac{E(u)}{u} (u-v)\Big)
\frac{dv}{v}.
$$
After the changes of variables $v=ut/3e^{\xi}$ and
$v=ut/E(u)$, respectively, this
becomes
$$
\align
& \gamma(E(u)) +\xi \log\Big(\frac{2}{3} \frac{E(u)}{E(u)-1}\Big)
+\int_{1}^{e^{\xi}} \log\Big(\frac{E(u)(1-t/(3e^{\xi}))}{E(u)-t}\Big)
\frac{dt}{t} \\
&\ge \gamma(E(u))
+\xi \Big(\log\Big(\frac{2}{3} \frac{E(u)}{E(u)-1}\Big) + \log \Big(
\frac{E(u)(1-1/(3e^{\xi}))}{E(u)-1}\Big)\Big),\\
\endalign
$$
since
$\log(\frac{E(u)(1-t/(3e^{\xi}))}{E(u)-t})$
is an increasing function of $t$, as
$3e^{\xi}>3>E(u)$.

Inserting the above bound for $I_2$ in (6.3) we deduce that
$$
\sigma(u)\ge
\rho(E(u)) + \frac{\xi}{2}\Big\{ (1-\xi)\Big( \log\Big(\frac{2}{3}
\frac{E(u)}{E(u)-1}\Big) + \log \Big(
\frac{E(u)(1-1/(3e^{\xi}))}{E(u)-1}\Big)\Big)-\gamma(E(u))\Big\},
$$
since
$\rho(x)=1-\log x+\gamma(x)/2$
in the range $2\le x\le 3$. Now, the quantity in
$\{\}$ is a decreasing function of
$E(u)$ (since each term is), and so is bounded
below by the value when substituting $2.6$ in for
$E(u)$, and this is positive for all
$\xi \le \log (1.6)$. It follows $\sigma(u)\ge \rho(E(u))$
and strict inequality holds unless $\xi =0$.  If $\xi=0$ then
$\sigma(u) = 1-I_1(u)+I_2(u)/2$ and if this equals $\rho(E(u))$
then one must have $I_2(u)=\gamma(E(u))$, and arguing
as above using Lemma 5.1, we see that this implies that $E(u_0)=1$.

\enddemo

\head 7. The case $E(u_0) \le E(u)-1$ and $E(u)>2.6$
\endhead

\noindent We call $u$ a ``{\sl champion for $\sigma$}''
if the absolute minimum of $\sigma(v)-\rho(E(v))$
in the interval $0\le v\le u$ is attained at $u$. Evidently
we need only establish Theorem 5 for champion $u$.

\proclaim{Proposition 7.1}  If $u$ is a champion for $\sigma$  and
$$
u\rho(E(u)) \le \int_0^u \chi(t) \rho(E(u-t)) dt , \tag{7.1}
$$
then $\sigma(u) \ge \rho(E(u))$.  Further if strict inequality holds in {\rm (7.1)}
then $\sigma(u)> \rho(E(u))$.
\endproclaim
\demo{Proof} Since $u$ is a champion for $\sigma$, we have
$$
\sigma(u)-\rho(E(u)) \le \sigma(v) - \rho(E(v))
$$
for all $0\le v\le u$.  Multiplying both sides by $\chi(u-v)$
and then integrating with respect to $v$ from $0$ to $u$, we obtain
$$
B(u) \Big(\sigma(u)-\rho(E(u))\Big)
\le u\sigma(u) - \int_0^u \chi(v) \rho(E(u-v)) dv
\le u\Big(\sigma(u) - \rho(E(u))\Big),
$$
by (2.1) and (7.1). The result follows as $B(u) \le u$.
\enddemo

We will complete the proof of Theorem 5 by showing that
(7.1) holds for $E(u)>2.6$ and $E(u_0)\le E(u)-1$ and
also determining when equality holds in (7.1).
We define
$$
I_1 = \int_{u_1}^u \rho(E(u-t)) \chi(t) dt,
\ \
I_2 = \int_{u_0}^{u_1} \rho(E(u-t)) \chi(t)dt,
\ \
\text{and}
\ \
I_3 = \int_{0}^{u_0} \rho(E(u-t)) \chi(t) dt.
$$
Put $V=u-B(u_0)$. Note that $B(u_0)\leq u_0$ and so $u\geq V\geq
u-u_0=u_1$.

 Since $\rho(E(u-t))$ is a non-decreasing function of $t$
we see by Lemma 5.1 that
$$
I_3
=\int_0^{u_0} \chi(t) \rho(E(u-t)) dt
\ge \int_{0}^{B(u_0)} \rho(E(u-t)) dt =
\int_{V}^{u} \rho(E(t)) dt.
$$
In the range $V\le t \le V E(u)/E(V)$
we have the bound $E(t) \le E(V) t/V$, and in
the range $VE(u)/E(V) \le t\le u$ we have the trivial bound $E(t) \le E(u)$.
Employing these bounds above we deduce that
$$
\align I_3 &\ge
\Big( u- \frac{VE(u)}{E(V)}\Big)\rho(E(u)) + \int_{V}^{{VE(u)/E(V)}}
\rho \Big( \frac{E(V)}{V} t \Big) dt \\
&= \Big( u- \frac{VE(u)}{E(V)}\Big)\rho(E(u))
+ \frac{V}{E(V)} \int_{E(V)}^{E(u)}
\rho(t) dt \\
&= u \rho(E(u)) - \frac{V}{E(V)} \int_{E(u)-1}^{E(V)} \rho(t) dt.\tag{7.2}
\\
\endalign
$$

Next
$$
\align
I_1 &\ge \rho(E(u_0)) (B(u) -B(u_1)) \ge
 \frac{\rho(E(u_0))}{\rho(E(u)-1)} (B(u)-B(u_1))  \rho(E(u)-1)   \\
&\ge \frac{V}{E(V)} \int_{E(u)-1}^{E(u_1)+\tau} \rho(t) dt, \tag{7.3a}\\
\endalign
$$
for $\tau$ satisfying
$$
(E(u_1)+\tau) -(E(u)-1) = \frac{E(V)}{V} \frac{\rho(E(u_0))} {\rho(E(u)-1)}
(B(u)-B(u_1)). \tag{7.3b}
$$
By Lemma 5.2 and (5.1) we have
$$
B(u) -B(u_1)
\ge \frac{uE(u_1)}{E(u)} - u_1 = \frac{u}{E(u)} (E(u_1) - E(u)+1) \ge
\frac{V}{E(V)} (E(u_1) -E(u)+1),
$$
and so $\tau\geq 0$ as $E(u_0)\leq E(u)-1$.

Finally note that
$$
I_2 \ge \rho(E(u_1)) (B(u_1)-B(u_0)) \ge \frac{V}{E(V)}
\int_{E(u_1)+\tau}^{E(u_1)+\tau+\tau'} \rho(t) dt, \tag{7.4}
$$
where $\tau'= \frac{E(V)}{V}(B(u_1)-B(u_0))$.

Combining the lower bounds given above for
$I_1, I_2$ and $I_3$, we see that (7.1)
follows provided $\tau+\tau'\geq E(V)-E(u_1)$.
Now let $C$ be a real number $\geq 1$ such that
$\rho(E(u_0))\geq C \rho(E(u)-1)$.
Define $\eta:=u_0 E(V)/V\geq 1$ by (5.1), and
$\lambda:=E(u)-1/E(u_0)$.
By (5.1) we have $V=u-B(u_0) \le u - u_0/E(u_0)
=u_0\lambda$, and so $\eta \lambda \geq \eta V/u_0=E(V)$.
Therefore if
$$
E(u_0) \ge \eta(1+C-C\eta) \lambda + 1+\eta -E(u)\tag{7.5}
$$
then, by Lemma 5.2,
$$
\align \frac{V}{E(V)} (\tau+\tau'+(E(u_1)-E(u)+1))
&= \frac{\rho(E(u_0))} {\rho(E(u)-1)} (B(u)-B(u_1)) + (B(u_1)-B(u_0))\\
&\geq C (B(u)-B(u_1)) + (B(u_1)-B(u_0))\\
&\ge C(B(u)-B(V)) + B(V)-B(u_0)   \\
&\ge C\Big(u\frac{E(V)}{E(u)}-V\Big)+ V\frac{E(u_0)}{E(V)} -u_0\\
&\ge \frac{V}{E(V)} (E(V)-E(u)+1), \\ \endalign
$$
so that $\tau+\tau'\geq E(V)-E(u_1)$ as desired.
Further if strict inequality holds in (7.5) then the inequality
in (7.1) is also strict.

If $C\ge 1+1/\lambda$ then the right side of (7.5) is
decreasing in $\eta\ge 1$, so
that it suffices to verify (7.5) at $\eta=1$.
This states that $E(u_0) \ge \lambda
+2-E(u) = 2-1/E(u_0)$ which always holds, further
the inequality is strict unless $E(u_0)=1$.
Consequently if $\rho(E(u_0))/\rho(E(u)-1)\ge E(u)/(E(u)-1)$ then
criterion (7.1) follows, since $\lambda \ge E(u)-1$, and further
(7.1) holds strictly unless $E(u_0)=1$.

If $C<1+1/\lambda$, then the right side of (7.5)
attains its maximum when $\eta =
(C+1+1/\lambda)/(2C)$, so that (7.5) holds if
$$
E(u_0) + \frac{(C+1)^2}{4CE(u_0)}
\ge E(u) \frac{(C-1)^2}{4C} + 1 + \frac{C+1}{2C}
+\frac{1}{4C\lambda}. \tag{7.6}
$$
Taking $C=1$ and noting that $\lambda\geq 2.6-1/E(u_0)$ we find
that strict inequality in (7.6) holds (and thus strict inequality in (7.1))
if $E(u_0)\ge 1.4341$. Hence we may assume that $E(u_0)<1.4341$.
If $E(u) \ge 2.802$ then $\rho(E(u_0))/\rho(E(u)-1) \ge
\rho(1.4341)/\rho(1.802)
>2.802/1.802\geq E(u)/(E(u)-1)$ so that (7.1) holds (and
equality there is possible only when $E(u_0)=1$).
Hence we may assume that $2.6\leq E(u) < 2.802$.

Now take $C=\rho(1.4341)/\rho(1.6)>1$ so that (7.6)
holds strictly (and thus (7.1) holds strictly) if $E(u_0)\ge
1.2383$. Hence we may assume that $E(u_0)<1.2383$. If $E(u) \ge 2.6635$ then
$\rho(E(u_0))/\rho(E(u)-1) \ge \rho(1.2383)/\rho(1.6635)
>2.6635/1.6635\geq E(u)/(E(u)-1)$ so that (7.1) holds (again
with equality only when $E(u_0)=1$).  Hence we may assume that
$2.6\leq E(u) < 2.6635$.

Now take $C=\rho(1.2383)/\rho(1.6)>1$ so that
(7.6) holds strictly (and thus (7.1) strictly) if $E(u_0)\ge
1.0648$. Hence we may assume that $E(u_0)<1.0648$. If $E(u) \ge 2.6$ then
$\rho(E(u_0))/\rho(E(u)-1) \ge \rho(1.0648)/\rho(1.6)
>2.6/1.6 \geq E(u)/(E(u)-1)$ so that (7.1) holds and with
equality possible only when $E(u_0)=1$.

\enddocument